\documentclass[12pt]{amsart}
\usepackage[utf8]{inputenc}

\usepackage[margin=1in]{geometry}
\usepackage[english]{babel}
\usepackage{color}
\usepackage{chngcntr}
\usepackage{amssymb, amsmath, amsthm, framed, mathabx, amsfonts, relsize, bbm, float, thmtools, graphicx, caption, mathtools}
\usepackage{parskip}
\usepackage{mathrsfs, fancyhdr, hyperref}
\usepackage{tikz, mathtools}
\usepackage{enumitem}
\usepackage{multirow, multicol, tabulary}

\newtheorem{theorem}{Theorem}[section]
\newtheorem*{theorem*}{Theorem}
\newtheorem{corollary}[theorem]{Corollary}
\newtheorem{lemma}[theorem]{Lemma}

\newtheorem{prop}[theorem]{Proposition}
\newtheorem{remark}[theorem]{Remark}
\newtheorem{case}{Case}
\newtheorem{definition}[theorem]{Definition}

\newcommand{\ba}{\begin{aligned}}
\newcommand{\ea}{\end{aligned}}
\newcommand{\bas}{\begingroup \addtolength{\jot}{.35em} \begin{aligned}}
\newcommand{\eas}{\end{aligned}\endgroup}
\newcommand{\bal} {\renewcommand*{\arraystretch}{1}\begin{array}{l}}
\renewcommand{\bar} {\renewcommand*{\arraystretch}{1}\begin{array}{r}}
\newcommand{\bp}{\begin{proof}}
\newcommand{\ep}{\end{proof}\vspace{2em}}
\newcommand{\bl}{\begin{lemma}}
\newcommand{\el}{\end{lemma}}
\newcommand{\bi}{\begin{itemize}}
\newcommand{\ei}{\end{itemize}}
\newcommand{\be}{\begin{enumerate}}
\newcommand{\bea}{\begin{enumerate}[label = \alph*.]}
\newcommand{\ee}{\end{enumerate}}

\newcommand{\eps}{\varepsilon}
\renewcommand{\phi}{\varphi}

\newcommand{\R}{\mathbb{R}}

\newcommand{\id}{\mathbbm{1}}

\newcommand{\dv}{\;\text{d}vol_g}

\newcommand{\calO}{\mathcal{O}}
\newcommand{\calW}{\mathcal{W}}
\newcommand{\calF}{\mathcal{F}}

\newcommand{\grad}{\nabla}

\renewcommand{\div}{\mathrm{div}}

\definecolor{green}{rgb}{0.0, 0.5, 0.0}

\DeclareMathOperator{\ric}{Ric}
\DeclareMathOperator{\tr}{tr}
\DeclareMathOperator{\hess}{Hess}
\DeclareMathOperator{\Lie}{\mathcal{L}}

\counterwithin*{case}{prop}

\title{Gradient Ambient Obstruction Solitons on Homogeneous Manifolds}
\author{Erin Griffin}
\address{215 Carnegie Building\\
Dept. of Math, Syracuse University\\
Syracuse, NY, 13244.}
\email{\href{mailto:egriff02@syr.edu}{egriff02@syr.edu} }
\urladdr{\url{https://sites.google.com/view/erin-griffin-math/}}
\thanks{The author was partially supported by grant NSF-\#1654034.}

\begin{document}

\maketitle

\begin{abstract}
\noindent We examine homogeneous solitons of the ambient obstruction flow and, in particular, prove that any compact ambient obstruction soliton with constant scalar curvature is trivial. Focusing on dimension 4, we show that any homogeneous gradient Bach soliton that is steady must be Bach flat, and that the only non-Bach-flat shrinking gradient solitons are product metrics on $\R^2\times S^2$ and $\R^2 \times H^2$. We also construct a non-Bach-flat expanding homogeneous gradient Bach soliton. We also establish a number of results for solitons to the geometric flow by a general tensor $q$. 
\end{abstract}

\section{Introduction}\label{intro}

The geometric flow for a general tensor $q(g)$, the $q$-flow, is a one parameter family of smooth metrics such that
\begin{equation} \label{qflow}
\begin{cases} \partial_t g= q
\\g(0) = h.
\end{cases}
\end{equation}
The resulting $q$-soliton equation is:
\begin{equation}\label{qs} 
\frac{1}{2} \Lie_X g = cg + \frac{1}{2} q \end{equation}
where $X$ is a vector field. Letting $X= \grad f$, a (normalized) gradient soliton has the form:
\begin{equation} \label{gqs} 
\hess{f} = cg + \frac{1}{2}q.
\end{equation}
The coefficient on $q$ may differ from other definitions found throughout the literature. We chose this coefficient to show that gradient solitons are self similar solutions to the $q$-flow (Theorem \ref{solutions}). It is easily shown that this definition aligns with definitions lacking the coefficient. We will use the terms expanding, steady, and shrinking to describe when $c<0$, $c=0$, and $c>0$ respectively. Moreover, a soliton is said to be stationary if it has constant potential function. 

The goal of our work is to generalize results for specific flows using the properties of $q$, then show that there are examples of these generalizations. One such general result is:

\begin{theorem} \label{q gen}
For a divergence-free, trace-free tensor $q$, any compact $q$-soliton is $q$-flat.
\end{theorem}

This theorem is a generalization of the well known result for Ricci solitons that any compact Ricci soliton with constant scalar curvature is Einstein (see \cite{PWRigidity}). By the second Bianchi identity, the Ricci tensor is divergence free if and only if the scalar curvature is constant. Moreover, for the Ricci tensor requiring constant scalar curvature is similar to the trace-free condition.

Two examples of divergence-free and trace free tensors are the Bach tensor, $B$, in dimension $n=4$, and ambient obstruction tensor, $\calO$, in even dimension $n\geq 4$. To aid our exposition, we provide detailed definitions of these tensors, their flows, and their solitons in Section \ref{background}. Specifically, Definition \ref{solitons} defines ambient obstruction and Bach solitons. In applying Theorem \ref{q gen} to the Bach tensor we get \cite[Theorem 3.2]{Ho}. Moreover, applying Theorem \ref{q gen} to the ambient obstruction tensor yields:

\begin{theorem} \label{AO gen}
Any compact ambient obstruction soliton with constant scalar curvature is $\calO$-flat. In particular, all compact homogeneous ambient obstruction solitons are $\calO$-flat.
\end{theorem}

The following is a generalization of \cite[Theorem 1.1] {PWRigidity} applied to the ambient obstruction tensor.

\begin{theorem} \label{intro ric}
A compact gradient ambient obstruction soliton with non-positive Ricci curvature must be stationary.
\end{theorem}

Our results are motivated by a recent theorem of Petersen and Wylie in \cite{PW} that implies that non-flat homogeneous gradient solitons of the $q$-flow, where $q$ is divergence free, are always products of the form $\R^k \times N^{n-k}$. In the case of Ricci solitons, this implies that $N$ is an Einstein manifold \cite[Theorem 1.1]{PWSymmetry}. We apply this theorem to the case of the Bach tensor and find that the possible metrics on $N$ are more complex. However, we obtain the following classification in the steady and shrinking cases.

\begin{theorem}\label{Bach gen}
Any homogeneous gradient Bach soliton that is steady must be Bach flat and the only non-Bach-flat shrinking solitons are product metrics on $\R^2 \times S^2$ and $\R^2 \times H^2$.
\end{theorem}

\begin{remark}
There are non-trivial homogeneous $4$-dimensional Bach flat metrics. For example, Einstein metrics and (anti)self-dual metrics are Bach flat. Moreover there is a classification of simply connected homogeneous Bach-flat $4$-manifolds. (See \cite{AGS} and \cite{BachFlat}.)
\end{remark}

\begin{remark}
There are non-Bach-flat expanding homogeneous gradient Bach solitons. We find one such soliton on $\R\times S^3$ with metric $g= g_0 \times g_{SU(2)}$. We show this is the only expanding soliton on a manifold of the form $\R\times N^3$ where $N^3$ is a unimodular Lie group.
\end{remark}

The paper is organized as follows. In Section \ref{background} we provide a brief background of the Bach and ambient obstruction tensors as well as a background on geometric flows in general. Next, in Section \ref{gen tensor} we establish a number of results for a general tensor $q$ and apply them to the ambient obstruction tensor. Then in Section \ref{Bach Flow} we begin to classify the gradient Bach tensors of homogeneous 4-manifolds. The results of this partial classification are summarized in Table \ref{table}.

\section{Background}\label{background}

In dimension 4 the Bach tensor is symmetric, divergence free, trace free, and conformally invariant of weight -2. That is, for a positive, smooth function $\rho$,  if $\tilde{g} = \rho^2 g$ then $\tilde{B} = \frac{1}{\rho^2} B$.  The Bach tensor is realized as the negative gradient of the conformally invariant functional given by:
\[ \calW(g) = \int_M |W_g|^2 \text{d}V_g \]
where $W_g$ is the Weyl tensor and $|W_g|^2 = g^{ip} g^{jq} g^{kr} g^{ls} W_{ijkl} W_{pqrs}$. Since this functional is only conformally invariant in dimension $n=4$, for $n\neq 4$ the Bach tensor is not conformally invariant either. Moreover, the Bach tensor is not divergence free for $n\neq 4$. For this reason we will only consider the Bach tensor for $n=4$. Though there is an explicit representation of the Bach tensor for arbitrary $n$, provided in \cite{CaoChen}, the Bach tensor for $n=4$ is given by:
\[ B_{ij} = \nabla^k \nabla^l W_{ikjl} + \frac{1}{2} R_{kl} W_{i\;\; j}^{\;k\;\, l}. \] 

To find a higher dimensional equivalent, for even $n$ we examine the gradient of the functional:
\[ \calF_Q^n (g) = \int_M Q(g) \,\text{d}V_g\]
where $Q(g)$ is Branson's $Q$-curvature described in \cite{Branson}. Though $Q$ lacks some of the conformal properties of $W$, the functionals $\calF_Q^n$ are conformally invariant for arbitrary even $n$. Moreover, Branson uses the Chern-Gauss-Bonnet theorem to show that, in dimension $n=4$, $\calF_Q^n$ is related to $\calW$ by the equation:
\[\calF_Q^4 = 8\pi^2 \chi(M) - \frac{1}{4} \calW,\]
thus they have the same critical metrics.

In \cite{FGAmbient}, Fefferman and Graham examine the gradient of $\calF_Q^n$ and introduce the resulting symmetric 2-tensor, the ambient obstruction tensor, which is noted $\calO$. This tensor can be characterized as the obstruction to an $n$-manifold having a formal power series of asymptotically hyperbolic Einstein metric in dimension $n+1$ \cite{BHexistence}. Explicitly, the ambient obstruction tensor is given by the equation:
\[ \calO_n = \frac{1}{(-2)^{\frac{n}{2}-2} \left(\frac{n}{2}-2 \right)! } \left( \Delta^{\frac{n}{2}-1} P - \frac{1}{2(n-1)} \Delta^{\frac{n}{2}-2} \nabla^2 S \right) + T_{n-1}\]
\[ P = \frac{1}{n-2} \left( \ric - \frac{1}{2(k-1)} S_g \right)\]
where $P$ is the Schouten tensor and $T_{n-1}$ a polynomial natural tensor of order $n-1$. The ambient obstruction tensor is only defined for even $n$. Like the Bach tensor in dimension $4$, the ambient obstruction tensor is symmetric, trace free, divergence free, and conformally invariant of weight $2-n$. The ambient obstruction tensor can be viewed as a family of even dimensional tensors, where the dimension 4 ambient obstruction tensor is the Bach tensor. (See \cite{BHexistence} and \cite{Lopez} for a more detailed background.)

In the last decade Bahuaud-Helliwell, Helliwell, and Lopez have studied flowing a metric by the ambient obstruction tensor. Bahuaud and Helliwell, in \cite[Theorem C]{BHexistence}, consider the flow given by:

\begin{equation} \label{AOF}
\begin{cases} 
\partial_t g = \calO_n +c_n (-1)^{\frac{n}{2}}\left(\Delta^{\frac{n}{2}-1} S\right)g \\
g(0)=h
\end{cases}\end{equation}
where $h$ is a smooth metric on a compact manifold of even dimension $n\geq 4$ and 
\[c_n = \frac{1}{2^{\frac{n}{2}-2}\left( \frac{n}{2} - 2\right)!(n-2) (n-1) }. \]
In \cite{BHexistence, BHuniqueness} Bahuaud and Helliwell show short time existence and uniqueness on compact manifolds for this flow. As Lopez explains in \cite{Lopez}, the scalar curvature term ``counteracts the invariance of $\calO$ under the action of the conformal group on the space of metrics on $M$.'' In \cite{Lopez}, Lopez finds pointwise smoothing estimates and uses them to find an obstruction to long-time existence and to prove a compactness theorem for the flow (\ref{AOF}). 
 
For $n=4$ we will call flow (\ref{AOF}) the Bach flow, which is given by: 
\[ \begin{cases} \partial_t g = B+ \frac{1}{12} \Delta S g 
\\ g(0) = h \end{cases} \]
Note that this is slightly different than the definition in \cite{Ho}. Since homogeneous manifolds have constant scalar curvature, the equations for the ambient obstruction flow and Bach flow on homogeneous manifolds are given by:
\begin{equation} \label{flow}
\begin{cases} \partial_t g= \calO_n
\\g(0) = h 
\end{cases} \quad\text{and} \quad
\begin{cases} \partial_t g= B
\\g(0) = h 
\end{cases}\end{equation}
respectively. Helliwell uses the latter equation in \cite{Helliwell} to study the Bach flow on homogeneous compact product manifolds of the form $S^1 \times K^3$.

The solitons of these flows are defined as follows.
\begin{definition}\label{solitons} An ambient obstruction soliton is a solution, $(M,g)$, to the equation:
\[ \frac{1}{2} \Lie_X g = cg + \frac{1}{2} \left( \calO_n +c_n (-1)^{\frac{n}{2}}\left(\Delta^{\frac{n}{2}-1} S\right) g\right)  \]
where $c_n$ is defined as above. In dimension $n=4$, the ambient obstruction soliton is the Bach soliton, given by:
\[ \frac{1}{2} \Lie_X g = cg + \frac{1}{2} \left( B +\frac{1}{12} \Delta S g\right). \]
These are called gradient if $X= \nabla f$, and the corresponding equations are
\[\hess{f} = cg + \frac{1}{2} \left( \calO_n +c_n (-1)^{\frac{n}{2}}\left(\Delta^{\frac{n}{2}-1} S\right)g\right) \quad \text{and}\]
\[\hess{f} = cg + \frac{1}{2}\left( B +\frac{1}{12} \Delta S g\right).\]
\end{definition}

For the reader who is less familiar with geometric flows, we now give a brief background will help motivate  gradient solitons.

The primary objective of a first course in differential equations is learning methods to solve differential equations explicitly. Soon thereafter, we see that solvable differential equations are relatively rare. To gain valuable insights about a differential equation, one might examine the fixed points of the flow, classify them as stable or unstable, and even construct a phase diagrams. Applying this idea to geometric flows, we know that self similar solutions to a geometric flow are solitons. As such they act as fixed points. So when examining a new flow it makes sense to try to find and analyze the solitons. To further limit these unknowns, one might choose to examine gradient solitons in particular. 

The choice to examine gradient solitons provides one with a more restrictive, more familiar environment. Historically, analyzing gradient solitons has provided a lot of insight into the Ricci flow. The work of Hamilton, Ivey, and Perelman combine to classify 3-dimensional shrinking gradient Ricci solitons.\cite{PWIntro}  Further, in \cite{Perelman}, Perelman observes that any compact Ricci soliton is a gradient Ricci soliton. Most notably, the study of Ricci solitons was imperative in Perelman's proof of the Poincar\'e Conjecture. The study of Ricci solitons has continued to prove a bountiful source of information and is still a very large area of research. It is reasonable to hope that the study of gradient solitons for other flows (specifically the Bach flow and ambient obstruction flow) would prove similarly fruitful in the understanding of the behavior of the flows and consequently the behaviors of the tensors themselves.

\section{Results for General Tensor}\label{gen tensor}

In this section, we prove a number of statements for a general trace free and/or divergence free tensor $q$. Applications of the theorem to the ambient obstruction tensor will follow in subsequent corollaries. For the sake of simplicity, full proofs of these corollaries have been omitted, but appropriate connections will be made.

Recall from Section \ref{background} that the ambient obstruction tensor, $\calO_n$ $n$ even, is trace free and divergence-free. However, the reader should note that the tensor affiliated with the general flow (\ref{AOF}) does not possess all of these properties.

One fact that proves useful in examining gradient solitons is the following proposition.

\begin{prop}  \label{ric grad f G}
Let $q$ be a symmetric two tensor and $(M,g,f)$ a gradient $q$-soliton (\ref{gqs}). The potential function, $f$, has the property that 
\[ \ric({\nabla f}) = \div{Q} - \frac{1}{2}\nabla(\tr {Q}) \] 
where $Q$ is the dual (1,1)-tensor of $q$ with respect to $g$.
\end{prop}

\bp
Consider a gradient soliton of the  $q$-flow, given by 
\[ \hess{f} = cg_{ij} + \frac{1}{2} q_{ij} \]
Type changing into (1,1) tensor
\[ \nabla \nabla f = cI + \frac{1}{2} Q \]
If we simply take the trace of each of the terms, we see that then  $\Delta f = cn + \frac{1}{2}\tr {Q}$.

Taking the divergence of each term in our soliton equation we see that:
\[ \ba
 \div Q &= \div (\nabla \nabla f)
\\&= \ric({\nabla f}) + \nabla( \Delta f) 
\\&= \ric({\nabla f}) + \nabla( cn + \frac{1}{2} \tr{Q})
\\&=\ric({\nabla f}) + \frac{1}{2}\nabla(\tr {Q})
\ea\]
Thus:
\[ \ric({\nabla f}) = \div{Q} - \frac{1}{2}\nabla(\tr {Q}) \] 
\ep

This theorem can be used to generalize \cite[Theorem 3.4]{Ho} as follows.
\begin{corollary} \label{ric zero}
For any constant trace, divergence free tensor $q$, the gradient solitons of its flow has that property that $\ric(\nabla f) = 0$ 
\end{corollary}

For the ambient obstruction flow on a non-homogeneous manifold, we see that a gradient soliton is given by:
\[ \hess{f} = cg + \frac{1}{2} \left(\calO_n +a_n \left(\Delta^{\frac{n}{2}-1} S\right)g \right) \]
where 
\[a_n = \frac{(-1)^{\frac{n}{2}}}{2^{\frac{n}{2}-2}\left( \frac{n}{2} - 2\right)!(n-2) (n-1)}.\] 
Note that $a_n$ simply combines constant terms in our original definition to help with notation. Examining this soliton, we get the following corollary. 

\begin{corollary} \label{AOF Ric}
A gradient ambient obstruction soliton with potential function $f$ satisfies $\ric(\nabla f) = a_n (1-n) \;\nabla \left( \Delta^{\frac{n}{2}-1} S\right) $. 
\end{corollary}

\bp
Consider a gradient ambient obstruction soliton with potential function $f$. Then $q=\calO_n +a_n \left(\Delta^{\frac{n}{2}-1} S\right)g $ and consequently
\[\ba \div{q} &= a_n d\left( \Delta^{\frac{n}{2}-1} S\right) 
 \\ \tr{q} &= n a_n \left( \Delta^{\frac{n}{2}-1} S\right) 
 \\ \nabla \tr{q} &= n a_n \nabla \left( \Delta^{\frac{n}{2}-1} S\right).
 \ea\]
Using Proposition \ref{ric grad f G}:
 \[ \ric({\nabla f}) = a_n (1-n) \nabla \left( \Delta^{\frac{n}{2}-1} S\right) \] 
\ep

\begin{remark} \label{aof ric zero}
For a gradient ambient obstruction soliton with constant scalar curvature (specifically for homogeneous manifolds) we see that $\Delta^{\frac{n}{2}-1} S = 0$, so $\ric({\nabla f}) = 0$.
\end{remark}

The following lemma appears to be well known, but we include the proof for completeness.

\begin{lemma} \label{Lie identity}
For any symmetric (0,2)-tensor field $\psi$ and vector field $\xi$:
\[ \langle \Lie_\xi g, \psi\rangle = 2\div (i_\xi \psi) - 2 (\div \psi)\xi\]
where $i_\xi \psi$ is a 1-form such that $i_\xi \psi(\cdot) = \psi(\xi, \cdot)$
\end{lemma}

\bp
Consider a symmetric (0,2)-tensor field $\psi$ and a vector field $\xi$. For a (0,2)-tensor $A$, we know that $A(x,y) = g(A(x), y)$, so:
\[ \langle A, B \rangle = \sum_i g( A(e_i), B(e_i)) = \sum_i A(e_i, B(e_i))\]
where $B$ is a (1,1)-tensor.

Consider the Lie derivative as our (0,2)-tensor, and $\psi$ a (1,1)-tensor. First, examining the type change, consider $\psi$ as a (0,2)-tensor:
\[ \psi(X, Y) = g(\psi(X), Y)  \Longrightarrow \psi(X, E_j) = g(\psi(X), E_j) \Longrightarrow \psi(X) = \sum_j g(\psi(X), E_j) E_j\] 
Next, we know that: 
\[ \div(\iota_\xi \psi) =  \sum_i (\nabla_{E_i} \iota_\xi \psi)(E_i) = \sum_i \nabla_{E_i} \psi\left(\xi, E_i\right) =\sum_i \nabla_{E_i} g\left(\psi(E_i), \xi\right)\]
\[ (\div \psi)(\xi) =\sum_i g(\xi, \nabla_{E_i} (\psi(E_i)) )\]
Then
\[\begin{aligned}
\langle \Lie_\xi g, \psi \rangle &= \sum_i \Lie_\xi g(E_i, \psi(E_i)) 
\\&= \sum_i g\left( \nabla_ {E_i} \xi, \psi(E_i)\right)+ \sum_i g\left( E_i, \nabla_{\psi(E_i)} \xi\right) 
\\&= \sum_i g\left( \nabla_ {E_i} \xi, g(\psi(E_i), E_j) E_j \right)+ \sum_i g\left( E_i, \nabla_{ g(\psi(E_i), E_j) E_j } \xi\right)  
\\&= \sum_i g( \psi(E_i), E_j) g(\nabla_{E_i} \xi, E_j) + \sum_i g( \psi(E_i), E_j) g(E_i, \nabla_{E_j}\xi) 
\\&= 2 g( \psi(E_i), E_j) g(\nabla_{E_i} \xi, E_j)
\\&=2(g(\nabla_{E_i} \xi , \psi(E_i))
\\&= 2\left[ \nabla_{E_i} g(\xi, \psi(E_i)) - g(X, \nabla_{E_i}(\psi(E_i))) \right]
\\&= 2 \div \iota_\xi \psi - 2(\div\psi)(\xi)
\end{aligned}\]
Thus, the identity holds.
\ep

We use this fact to prove the following lemma for compact solitons of a general $q$-flow. 

\begin{lemma} \label{Lie div}
Let $(M, g, X)$ be an $n$-dimensional compact soliton to the $q$-flow, (\ref{qs}). Then:
\bea
\item $\mathlarger{ \int_M ||\Lie_X g||^2 \dv= -2\int_M \div(q)(X)\dv }$.
\item If $q$ is divergence free, then $X$ is Killing. 
\item If $q$ is divergence free and trace free, then $(M, g_{ij})$ must be $q$-flat. 
\ee
\end{lemma}

\bp
\bea 
\item Consider the gradient $q$-soliton, $\frac{1}{2} \Lie_X g = cg + \frac{1}{2} q$. We know that for any vector field $\xi$ on $M$
\[ \langle \Lie_\xi g, \psi\rangle = 2 \div(i_\xi \psi) - 2 (\div \psi) (\xi) \]
where $i_\xi\psi(\cdot) = \psi(\xi, \cdot)$.

Note that the soliton can be written as $q = \Lie_X - 2cg$. Examining the divergence of this equation:
\[ \div q_{ij} = \div(\Lie_X g) - 2c \div( g_{ij})=  \div(\Lie_X g)\]
Using Lemma \ref{Lie identity}, we see that letting $\psi = \Lie_X g$ and $\xi = X$:
\[ \langle \Lie_X g, \Lie_X g \rangle = ||\Lie_X g||^2 = 2 \div (i_X \Lie_x g) - 2 \div(\Lie_X g)(X) = 2 \div (i_X \Lie_x g) - 2 \div(q)(X)\]

Integrating over $M$ we see that since $M$ is compact and has no boundary:
\[  \int_M ||\Lie_X g||^2 \dv= 2\int_M   \div (i_X \Lie_x g) \dv - 2\int_M  \div(q)(X)\dv = -2\int_M \div(q)(X)) \dv \]

\item If $q$ is divergence free part (a) shows that $\int_M ||\Lie_X g||^2 \dv = 0$. Thus, $\Lie_X g = 0$ and consequently $X$ is Killing.  

\item Suppose that $q$ is divergence free and trace free. From (b), this means that $q_{ij} = cg_{ij}$. Taking the trace of both sides we see that $0 = nc$ and thus $c=0$. Thus $q_{ij} = 0$ and subsequently $(M, g_{ij})$ is $q$-flat.
\ee \ep

\begin{corollary}  \label{aof compact}
Let $(M, g, X)$ be an $n$-dimensional compact soliton to the ambient obstruction flow with constant scalar curvature. Then $\int_M ||\Lie_X g||^2 \dv= 0$, $X$ is Killing, and $M$ is $\calO$-flat. 
\end{corollary}

\bp
Since $M$ has constant scalar curvature we know that the flow is given by (\ref{flow}). Thus, we consider $q= \calO_n$. Since $\calO$ is divergence free and trace free, the conclusion follows directly from Lemma \ref{Lie div}
\ep

In particular, Corollary \ref{aof compact} shows that any homogeneous compact ambient obstruction soliton is $\calO$-flat. In the non-homogeneous gradient case we have the following inequality.

\begin{theorem} \label{ric compact}
For any compact gradient ambient obstruction soliton $(M,g,f)$ 
\[\int_M \ric(\nabla f, \nabla f) \dv \geq 0,\]
where the integral is zero if and only if $f$ is constant.
\end{theorem}

\bp
Consider an $n$-dimensional compact gradient ambient obstruction soliton, $(M, g, f)$. Applying Lemma \ref{Lie div}, let $q= \calO$ and let $X = \nabla f$. From Corollary \ref{AOF Ric}: 
\[\div{Q} = a_n \nabla \left( \Delta^{\frac{n}{2}-1} S\right) = \frac{a_n}{1-n} (1-n) \nabla \left(\Delta^{\frac{n}{2}-1} S\right)= \frac{1}{1-n} \ric({\grad f}).\]
By Lemma \ref{Lie div}:
\[0 \leq   \int_M ||\Lie_{\grad{f}} g||^2 \dv= -2\int_M \div(q)(\grad{f})) \dv = \frac{2}{n-1} \int_M \ric\left(\grad{f}, \grad{f}\right) \dv.\]
Thus $\int_M \ric(\nabla f, \nabla f) \dv \geq 0$.

Suppose $\int_M \ric(\nabla f, \nabla f) \dv =0$. 

Since 
\[\int_M ||\Lie_{\grad{f}} g||^2 \dv = \frac{2}{n-1} \int_M \ric\left(\grad{f}, \grad{f}\right) \dv,\]
if the right hand side is zero then $\Lie_{\nabla f}(g)=0$ and consequently $\hess f = 0$. Since $M$ is compact this implies that $f$ is constant. If $f$ is constant $\nabla f = 0$ then clearly $\ric(\nabla f) = 0$. Therefore, the integral is zero if and only if $f$ is constant
\ep

\begin{remark}
Note that a soliton is defined to be stationary if $f$ is constant. Thus Theorem \ref{ric compact} implies Theorem \ref{intro ric}.
\end{remark}

We note that in general, stationary gradient ambient obstruction solitons are characterized by the following proposition. 

\begin{prop} \label{stationary}
If $(M,g,f)$ is a stationary gradient ambient obstruction soliton, then $(M,g)$ is $\calO$-flat. If $(M,g)$ is also compact then $S$ is constant.
\end{prop}

\bp

Consider a stationary gradient ambient obstruction soliton, $(M,g,f)$. Since the soliton is stationary, $f$ is constant. Consequently $\hess f = 0$ and thus  $q = -2cg$. Since $q= \calO_n+ a_n\left(\Delta^{\frac{n}{2}-1} S\right)$,
\[\calO_n = \left( -a_n \left(\Delta^{\frac{n}{2}-1} S\right) -2c\right) g.\]
Taking the trace of both sides:
\[ 0 = n\left( -a_n \left(\Delta^{\frac{n}{2}-1} S\right) -2c\right)\]
Thus
\[0 = -a_n \left(\Delta^{\frac{n}{2}-1} S\right) -2c\]
This forces $\calO_n =0$, so that soliton is $\calO$-flat. Furthermore:
\[ \Delta^{\frac{n}{2}-1} S = \frac{2c}{a_n} \]
is constant. If $M$ is compact, this implies that $S$ is constant. 
\ep

\begin{remark}
The converse of Proposition \ref{stationary} is true in the compact case. That is, a compact gradient ambient obstruction soliton that is $\calO$-flat and has constant scalar curvature is stationary. Constant scalar curvature and $\calO$-flat imply that $\hess f = cg$. Compactness forces the manifold to have a maximum and minimum so $\hess f = 0$. Appealing once more to compactness, this forces $f$ to be constant and our soliton to be stationary. 
\end{remark}

Though the following lemma is not necessary when studying ambient obstruction solitons (this was taken care of in Corollary \ref{aof compact}), it does give another criteria for when a $q$-soliton is stationary.

\begin{prop} \label{compact}
For a trace free tensor $q$, any compact gradient soliton to the $q$-flow must be $q$-flat.  
\end{prop}
\bp
Generalizing from \cite{Ho}, consider a gradient $q$-soliton (\ref{gqs}). By assumption $\tr(q) = 0$, so taking the trace of both sides yields $\Delta f =cn$. Integrating over $M$:
\[ 0 = \int_M c n - \Delta f \dv = cn \; Vol(M,g) \]
Thus $c = 0$. Further, $\Delta f = 0 n+ 0$ so $\Delta f = 0$, that is, $f$ is harmonic. Since $M$ is compact, $f$ must be constant.\\
 
Therefore $q_{ij} = 2 \hess{f}-2c g_{ij}  = 0$, so any compact gradient soliton is $q$-flat.
\ep

Proceeding, we will show that for a general tensor $q$ with certain scaling properties that a gradient $q$-soliton is a self similar solution to the $q$-flow. This observation appears to be made first by Lauret \cite{Lauret2}. To do so we will follow the proof from \cite[Chapter 4]{HRF} which shows that gradient Ricci solitons are self-similar solutions to the Ricci flow. Following our proof, we will apply the theorem to the ambient obstruction flow in both the homogeneous and non-homogeneous cases. In \cite{Lauret1} and \cite{Lauret2}, Lauret shows that the following theorem is true for general, non-gradient solitons and can be made into an if and only if statement. I have chosen to focus on the case of gradient solitons. Our goal in including the following proof is to motivate our choice to modify the equation for a soliton by including a factor of $\frac{1}{2}$ and to show a more explicit proof of this theorem. 

\begin{theorem}\label{solutions}
Consider any tensor $q$ with the property that when the metric is scaled by a constant $\lambda \in \R$:
\[ \tilde{g} = \lambda g \quad\Longrightarrow\quad \tilde{q} = \lambda^{\frac{w}{2}} q.\]

Consider a complete gradient $q$ soliton $(M^n, h, f_0, c)$, that is:
\[ \hess_h f_0 = c h + \frac{1}{2} q(h).\]
There exists an $\eps >0$ such that for all $t\in (-\eps, \eps)$ there is a solution $g_t$ of the $q$ flow with $g_0 = h$, diffeomorphisms $\phi_t$ with $\phi_0 = \id_{M^n}$, and functions $f(t) = f_t$ with $f(0)=f_0$, such that:

\begin{enumerate}
\item $\tau$ is scales the metric according to the function:
\[ \tau_t := \begin{cases} 
  e^{1-2ct} & w=2 \\
 \left(1- 2c\left(1-\frac{w}{2}\right)t\right)^{\frac{1}{1-\frac{w}{2}}} & w\neq 2,
\end{cases}\]

\item The vector field $X_t:= \tau_t^{\frac{w}{2}-1} \grad_h f_0$ exists,

\item $\phi_t: M^n \to M^n$ is the 1-parameter family of diffeomorphisms generated by $X_t$. So:
\[ \frac{\partial}{\partial t} \phi_t(x) = \tau_t^{\frac{w}{2}-1} \left( \grad_h f_0 \right) (\phi_t(x)),\]

\item $g_t$ is the pull back by $\phi_t$ of $h$ up to the scale factor $\tau_t$:
\[ g_t = \tau_t \phi_t^* h,\]

\item $f_t$ is the pull back by $\phi_t$ of $f_0$:
\[ f_t = f_0 \circ \phi_t = \phi_t ^*(f_0) .\]
\end{enumerate}
Moreover 
\[\hess_{g_t} f_t  = \frac{c}{\tau_t} g_t + \frac{1}{2} (q(g_t)) \]
or equivalently
\[ q(g_t) = -\frac{2c}{\tau_t}g_t + 2\hess_{g_t} f_t\]
and 
\[ \frac{\partial f}{\partial t} (t) = \tau^{\frac{w}{2}} \left| \grad_{g_t} f_t \right|_{g_t}^2.\] 
\end{theorem}

\bp
Construct a 1-parameter family of diffeomorphisms $\phi_t: M^n \to M^n$ generated by vector field $X_t = \tau^{\frac{w}{2}-1} \grad_h f_0$ defined for all $t$ such that $t\in (-\eps, \eps)$. Define $f_t = f_0 \circ \phi_t$ and $g_t = \tau_t \phi_t ^*h$.
\[ \left. \frac{\partial}{\partial t} \right|_{t=t_0} g_t  = \left. \frac{\partial}{\partial t} \right|_{t=t_0} (\tau_t \phi_t ^*h )  = \left(\frac{\partial}{\partial t} \tau_t \right) \phi_{t_0} ^*h   + \tau_{t_0} \frac{\partial}{\partial t} \big\vert_{t=t_0} \phi_t ^*h \] 
 
Using Remark 1.24 from \cite{HRF} we are able to assess the derivative of the pullback:
\[ \tau_{t_0} \left. \frac{\partial}{\partial t} \right|_{t=t_0} \phi_t ^*h =  \tau_{t_0} \Lie_{Y(t)} \left(\phi_{t_0}^* h\right) = \Lie_{Y(t)} \left(\tau_{t_0} \phi_{t_0}^* h\right)\]
 where 
 \[Y(t) := \left. \frac{\partial}{\partial t} \right|_{t=t_0}  \left( \phi_{t_0}^{-1} \circ \phi_t \right) =  (\phi_{t_0}^{-1})_*\left. \frac{\partial}{\partial t} \right|_{t=t_0} \phi_t. \]
 Note that for $\hat{g} = \lambda g$:
  \[ g(\nabla_g f, X) = df(X) = \tilde{g}(\nabla_{\hat{g}} f, X) = \lambda g(\nabla_{\hat{g}} f, X).\]
  So $\frac{1}{\lambda}\nabla_g f = \nabla_{\hat{g}} f$. Therefore:
  \[ \nabla_{g_{t_0}} f_{t_0} = \nabla_{\tau_{t_0} \phi_{t_0}^* h} f_{t_0} = \frac{1}{\tau_{t_0}}  \nabla_{ \phi_{t_0}^* h} f_{t_0} = \frac{1}{\tau_{t_0}}  \nabla_{ \phi_{t_0}^* h} \phi_{t_0}^* f_0= \frac{1}{\tau_{t_0}} \phi_{t_0}^*( \nabla_h f_0) = \phi_{t_0}^*\left( \frac{1}{\tau_{t_0}}\nabla_h f_0 \right). \]
Thus
 \[ \left. \frac{\partial}{\partial t} \right|_{t=t_0} \phi_t  = \tau_{t_0}^{\frac{w}{2}-1} \grad_{h} f_0 = \tau_{t_0}^{\frac{w}{2}}  \left(\frac{1}{\tau_{t_0}} \grad_{h} f_0 \right) =\tau_{t_0}^{\frac{w}{2}}  \left( (\phi_{t_0})_* \left(\grad_{g_{t_0}} f_{t_0} \right)\right). \]
Using this, we are able to evaluate the desired derivative and find one term of our initial sum:
 \[  \tau_{t_0} \left. \frac{\partial}{\partial t} \right|_{t=t_0} \phi_t ^*h = \tau_{t_0} \Lie_{Y(t)} \left(\phi_{t_0}^* h\right) = \Lie_{\tau_{t_0}^{\frac{w}{2}} \grad_{g_{t_0}}f_{t_0}} \left( \tau_{t_0} \phi_{t_0}^* h\right) = \tau_{t_0}^{\frac{w}{2}} \Lie_{\grad_{g_{t_0}}f_{t_0}}  g_{t_0} \]
 
To evaluate the derivative of $\tau$ we must consider each case.
\begin{case} For $w=2$ define $\tau_t = e^{1-2ct}$. Then:
 
 \[ \ba \left(\frac{\partial}{\partial t} \tau_t \right) \phi_{t_0} ^*h  &= -2c \tau \phi_{t_0} ^*h  
 \\ &= -2c g(t_0)  
 \ea\]
 \end{case}
 
 \begin{case} For $w\neq 2$ define $\tau_t = \left(1- 2c\left(1-\frac{w}{2}\right)t\right)^{\frac{1}{1-\frac{w}{2}}}$. We can compute the following:
 \[\ba \left(\frac{\partial}{\partial t} \tau_t \right) \phi_{t_0} ^*h &=  \frac{1}{1-\frac{w}{2}}  \left(1- 2c\left(1-\frac{w}{2}\right)t_0\right)^{\frac{1}{1-\frac{w}{2}}- 1} \left(-2c\left(1-\frac{w}{2}\right) \right) \left(\phi_{t_0} ^*h\right)  
    \\&= -2c \left(1- 2c\left(1-\frac{w}{2}\right)t_0\right)^{\frac{w/2}{1-\frac{w}{2}}} \left(\phi_{t_0} ^*h\right)
 \\ &= -2c\tau_{t_0}^{\frac{w}{2}} \left(\frac{\tau_{t_0}\phi_{t_0} ^*h}{\tau_{t_0}} \right) 
 \\ &= -2c\tau_{t_0}^{\frac{w}{2} -1} g(t_0) 
 \ea\]
\end{case}

Thus we see that for any $w$, 
\[ \left(\frac{\partial}{\partial t} \tau_t \right) \phi_{t_0} ^*h =-2c\tau_{t_0}^{\frac{w}{2} -1} g(t_0)  \]

Returning to our original derivative, we see that for general $t$:
\[\ba \frac{\partial}{\partial t} g_t =& -2c\tau_{t_0}^{\frac{w}{2} -1}  g_t + \tau_{t_0}^{\frac{w}{2}} \Lie_{\grad_{g_t} f_t} g(t) 
\\=& \tau_{t_0}^{\frac{w}{2}} \left( \frac{-2c}{\tau_t} g(t) + 2 \nabla^{g_t}\nabla^{g_t} f_t \right)\ea\]
Applying \cite{HRF} Exercise 1.23 to $q$ we see:
\[ \ba q(g_t) &= q ( \tau_t \phi_t^*h)
\\& = \tau_t^{\frac{w}{2}} \phi_t^* (q(h))
\\& =\tau_t^{\frac{w}{2}} \phi_t^*\left(-2c h + 2 \hess_h f_0 \right) 
\\& = \tau_t^{\frac{w}{2}} \phi_t^*\left(-2c h +  \Lie_{\grad_{h} f_0} h \right)
\\&= \tau_t^{\frac{w}{2}} \left(\frac{-2c}{\tau_t} g_t + \Lie_{\grad_{g_t} f_t} g(t) \right)
\\&= \tau_t^{\frac{w}{2}} \left(\frac{-2c}{\tau_t} g_t + 2\hess_{g_t} f_t\right)
\\&= \frac{\partial}{\partial t} g_t
\ea\]
Hence, there exists a solution $g_t$ to the flow with the desired properties. 

Looking at the derivative of the potential function we see that:
\[\ba \frac{\partial f_t(x)}{\partial t}  &= \frac{\partial}{\partial t} f_0(\phi_t(x))
\\&= \lim_{\eta\to 0} \frac{ f_0(\phi_{t+\eta}(x)) - f_0(\phi_t(x))}{\eta} 
\\&= h\left( \nabla_{h} f_0,\;  \frac{\partial}{\partial t} \phi_t \right)
\\&= h\left( \nabla_{h} f_0, \;  \tau^{\frac{w}{2}-1}\nabla_{h} f_0(\phi_t(x)) \right)
\\&= \tau^{\frac{w}{2}-1} h\left( \nabla_{h} f_t, \nabla_{h} f_t(x) \right)
\\&= \tau^{\frac{w}{2}-1}\; \frac{1}{\tau} g_t\left(\tau \nabla_{g_t} f_t, \; \tau \nabla_{g_t} f_t(x) \right)
\\&= \tau^{\frac{w}{2}} \left|\grad_{g_t} f_t \right|_{g_t}^2 
\ea\]

\ep

\begin{remark}
If the vector field $X_t=\tau_t^{\frac{w}{2}-1} \grad_h f_0$ is complete then the flow exists for all $t$ such that $\tau_t>0$. 
\end{remark}

\begin{remark}
One such tensor $q$ with the necessary weighting property is a conformally invariant tensor of weight $w$. That is, a tensor $T$ such that for $\tilde{g} = \rho^2 g$, then $\tilde{T} = \rho^w T$ for a smooth positive function $\rho$. 
\end{remark}

\begin{corollary} \label{aof solutions}
The gradient solitons of the ambient obstruction flow are self similar solutions to the ambient obstruction flow.
\end{corollary}
\bp
Consider the tensor provided by the ambient obstruction flow:
\[\calO_n + c_n(-1)^{\frac{n}{2}} \left( \Delta^{\frac{n}{2}-1} S\right)g.\]

We know that the ambient obstruction tensor is of conformal weight $2-n$, and is consequently a tensor $q$ described by Theorem \ref{solutions}. In the homogeneous case, or more generally the constant scalar curvature case, we are able to directly apply the theorem.

To examine the non-homogeneous case we must also investigate the scaling properties of the scalar curvature term. A simple calculation shows that for $\tilde{g} = \lambda^2 g$:
\[ \tilde{\Delta} \tilde{S} \tilde{g}= \frac{1}{\lambda ^2} \Delta S g .\]
Using induction one can show that this generalizes to:
\[\tilde{\Delta}^k \tilde{S} \tilde{g}  = \frac{1}{\lambda^{2k}} \Delta^k S g\]
Thus for $k = \frac{n}{2}-1$ 
\[\tilde{\Delta}^{\frac{n}{2}-1} \tilde{S} \tilde{g}  = \frac{1}{\lambda^{n-2}} \Delta^k S g = \lambda^{2-n} \Delta^k S g.\]
That is, the scalar curvature term is scaled by a factor of $2-n$ and consequently has the same scaling properties as the ambient obstruction tensor.

Applying Theorem \ref{solutions} with $w= 2-n$, we see that this implies that with the appropriate choice of $\tau$ and $\phi$ a gradient ambient obstruction soliton is a self-similar solution to the ambient obstruction flow.
\ep

As Lauret shows, Corollary \ref{aof solutions} is also true for non-gradient solitons. Turning our attention to noncompact, homogeneous solitons we consider recent theorem of Petersen and Wylie \cite{PW}. This theorem is a key part of understanding homogeneous gradient Bach solitons as we see in Section \ref{Bach Flow}. 

\begin{theorem} [Petersen-Wylie] \label{PW Theorem} 
Let $(M,g)$ be a homogeneous manifold and  $\hat{q}$ an isometry invariant symmetric two-tensor which is divergence free. If there is a non-constant function such that $\mathrm{Hess} f = \hat{q}$ then $(M,g)$ is a product metric $N\times\mathbb{R}^{k}$ and $f$ is a function on the Euclidean factor. 
\end{theorem}

For a divergence free tensor $q$, we apply this theorem to homogeneous gradient $q$ solitons by simply letting $\hat{q} = cg+ \frac{1}{2} q$. Then $\hat{q}$ is the sum of isometry invariant symmetric two-tensors that are divergence free and is itself such a tensor. Applying this theorem to homogeneous manifolds, we are able limit the ambient obstruction flow to the flow given by (\ref{flow}). Since $\calO$ is a divergence free, isometry invariant, symmetric two-tensor, we can let $q= \calO_n$ resulting in the following corollary. 

\begin{corollary}
If $(M,g)$ is a homogeneous gradient ambient obstruction soliton, then either $M$ is stationary or it splits as a product $\mathbb{R}^k\times N$ and $f$ is a function on the Euclidean factor. 
\end{corollary}

This theorem informs our approach to classifying homogeneous gradient Bach solitons in the next section.

\section{Gradient Bach Solitons}\label{Bach Flow}

In order to examine and classify the gradient solitons of the Bach flow on homogeneous 4-manifolds, we consider the four configurations of homogeneous 4-manifolds that are found by ``pulling off copies of $\R$''. More explicitly, by Theorem \ref{PW Theorem}, the solitons will be of the form $\R^4$, $\R^3 \times N^1$, $\R^2 \times N^2$, $\R \times N^3$, or $N^4$ (where $N^k$ is necessarily homogeneous). The first and last case we will call non-split manifolds, the others may be called the $3\times 1$, $2\times 2$, and $1\times 3$ cases respectively. For each of these cases (and for the remainder of the paper) it will be assumed that the product manifolds are equipped with the appropriate product metric $g=g_0 \times g_N$. Table 1 summarizes our findings regarding each type and thus proves the general theorem stated in the introduction. Prior to doing so, we set up the conventions used throughout this section. 

From (\ref{gqs}) we know that for homogeneous manifolds the equation for a gradient Bach soliton is given by: 
\[\hess f = cg + \frac{1}{2} B\]
and can be represented in coordinates as:
 \[ \nabla_i \nabla_j f = cg_{ij} + \frac{1}{2} B_{ij}. \]

In order to make the following proofs more clear, we will consider how the above equation can be given by matrices. In order to do this we will establish conventions that will hold for the remainder of the section unless otherwise noted. We will always choose a basis so both the metric and the Bach tensor are diagonal. (This is always possible, per the spectral theorem.) Since the metric and the Bach tensor are diagonal, $\hess f$ must also be diagonal so $\nabla_i \nabla_j f=0$ for $i\neq j$. One very important statement in Theorem \ref{PW Theorem} is that the potential function depends on only the Euclidean factor of the product manifold. Let $\nabla_i\nabla_i f = f_{ii}$. Thus, in general we see that the gradient Bach solitons can be represented by the following equality:
\[ \begin{bmatrix} f_{00}&0&0&0 \\ 0&f_{11}&0&0 \\ 0&0&f_{22}&0 \\ 0&0&0&f_{33}\end{bmatrix} = c \begin{bmatrix}g_{00}&0&0&0 \\ 0& g_{11} &0&0 \\ 0&0& g_{22} &0 \\ 0&0&0&g_{33} \end{bmatrix}    + \frac{1}{2} \begin{bmatrix}B_{00}&0&0&0 \\ 0& B_{11} &0&0 \\ 0&0& B_{22} &0 \\ 0&0&0&B_{33} \end{bmatrix}. \]

Recall from the introduction the generalization stated as Theorem \ref{Bach gen}. To prove this theorem we will simply examine each type of manifold and assess the solitons. The following table will summarize this investigation with one notable exception: in the $\R \times N^3$ case we are able to prove that non-Bach-flat gradient solitons must be expanding.

\begin{center}
\renewcommand{\arraystretch}{1.25}
\begin{tabulary}{1.07\textwidth}{|c|c|c|c|L|} 
\hline Split & Manifold &Type of Soliton &Permissible Metrics  & Potential Function\\\hline 
\multirow{3}{*}{$N^4$} 
    & $\R^4$ & Gaussian & Bach flat (any) &  $f(x,y,z,w) = c(x^2+y^2+z^2+w^2)+ax+by+dz+hw+k$ \\\cline{2-5} 
    & $N^4$ & Stationary &Bach flat & $f(x,y,z,w) = k$  \\\hline 

$\R^3 \times N^1$ & \;&  Steady & Bach flat (any)& $f(x,y,z) = ax+by+dz+k$   \\\hline 
\multirow{3}{*}{$\R^2 \times N^2$} 
    & $\R^2 \times \R^2$ & Steady & Bach flat (any)& $f(x,y) = ax+by+d$\\\cline{2-5} 
    & $\R^2 \times S^2$ & Shrinking & See \cite{Ho}& $f(x,y) = c(x^2+y^2)+ ax+by+dz+k$  \\\cline{2-5} 
    & $\R^2 \times H^2$ & Shrinking & See \cite{Ho}& $f(x,y) = c(x^2+y^2)+ax+by+dz+k$  \\\hline 
    
\multirow{10}{*}{$\R \times N^3$} & $\R \times \R^3$ & Steady & Bach flat (any) & $f(x) = ax+b$\\\cline{2-5} 
                     & $\R \times Nil$ & --- &None&  ---\\\cline{2-5} 
                     & $\R \times Solv$  & --- &  None & ---\\\cline{2-5} 
                     & $\R \times \widehat{SL}(2, \R)$ & --- & None&  --- \\\cline{2-5} 
                     & $\R \times (\R \times H^2)$ & --- & None &  ---\\\cline{2-5} 
                     & $\R \times (\R \times S^2)$ &---  & None&  --- \\\cline{2-5} 
                     & $\R \times E(2)$ &Steady  &  Bach flat ($g_{11}= g_{22}$)& $f(x) = ax+b$ \\\cline{2-5} 
                     & $\R \times H^3$ & Steady &  Bach flat & $f(x) = ax+b$\\\cline{2-5} 
                     & \multirow{2}{*} {$\R \times S^3$} &Steady  &  Bach flat ($g_{11} = g_{22}=g_{33}$)& $f(x) = ax+b$ \\\cline{3-5} 
                      & &Expanding&   $g_{11} = g_{22}=4g_{33}$& $f(x) = 2cx^2+ ax+b$ \\\hline 
\end{tabulary}
\captionof{table}{ Summary of Results }
\label{table}
\end{center}

\subsection{Non-split Manifolds} 

\begin{theorem} $(\R^4, g_0)$ is a Gaussian soliton. 
\end{theorem}

\bp 
We know from the equation for the Bach tensor that $(\R^4, g_0)$ is Bach flat, that is, $B_{ij} =0$ for all $i,j= 0,1,2,3$, so $\hess f = c g$. By Theorem \ref{PW Theorem}, $f$ is a function on $\R^4$. Thus for any orthonormal basis, $\R^4$ is a gradient Bach soliton with potential function 
\[f(x,y,z,w)= \frac{1}{2} c (x^2 + y^2+z^2+w^2)+ ax+by+dz+hw+k\]
for $a,b,d,h,k \in \R$. 

Since there are no restrictions on $c$, we see that this is the Gaussian soliton.
\ep

\begin{prop} 
Consider a non-split, homogeneous $4$-manifold $N^4 \neq \R^4$ with metric $g_N$. Then $N^4$ is a gradient Bach soliton if and only if it is Bach flat.
\end{prop}

\bp Consider a non-split, homogeneous $4$-manifold $N^4$ with metric $g_N$. By the converse of Theorem \ref{PW Theorem}, since $N^4$ is not a product manifold, it must have constant potential function and is therefore stationary. Since the potential function is constant, $\hess f =0$. Consequently, any soliton has the form $-\frac{1}{2}B= cg$. Taking the trace of each side we see that
\[ 0 = -\frac{1}{2} \tr{B} = \tr{cg} = 4c \]
and so it is necessarily true that $c=0$ and the soliton is steady.

Since $c=0$ always, $B=0$ always and thus the manifold must be Bach flat. 
\ep

\subsection{Manifolds of the form $\R^3 \times N^1$} 

\begin{remark}\label{R3N1}
For a manifold of the form $\R^3 \times N^1$ with metric $g= g_0 \times g_N$, we know that $N^1 = \R^1 \text{ or } S^1$. Thus any manifold of this form is flat and consequently Bach flat.
\end{remark}
 
\begin{prop}
Homogeneous manifolds of the form $\R^3 \times N^1$ with metric $g= g_0 \times g_N$ are steady gradient Bach solitons with linear potential functions. 
\end{prop}
\bp 
Consider a homogeneous manifold of the form $\R^3 \times N^1$ with metric $g= g_0 \times g_N$. We know from  Remark \ref{R3N1} that any manifold of this form is Bach flat. So for any gradient Bach soliton $\hess f = cg$. By Theorem \ref{PW Theorem}  we know that $f(x,y,z): \R^3 \to \R$. So $\nabla_3\nabla_3 f = 0 = cg_{33}$.  Since the metric is positive definite, $c=0$. Therefore, the gradient Bach solitons are steady. 

Consequently $\hess f = 0$,  so $f_{xx} = f_{yy}= f_{zz} = 0$. Thus $f(x,y,z) = ax+by+cz+d$. 
\ep

\subsection{Manifolds of the form $\R^2\times N^2$} 
In his 2018 paper, \cite{Ho}, Ho finds homogeneous gradient solitons of the form $\R^2 \times N^2$. Ho proves that both $\R^2 \times S^2$ and $\R^2 \times H^2$ is a nontrivial soliton of the form:
\[\hess f=  B +\frac{1}{12}g  \]
for any function $f$ of the form $f(x,y) = \frac{1}{12}(x^2 + y^2) + k$. Note the difference between Ho's definition of a gradient Bach soliton and that of this paper. Ho has chosen to place the metric term on the right hand side of the equation switching the conventions of shrinking/ expanding. We will prove that Ho's examples are the only examples of this type. 

\begin{theorem} \label{2x2 shrinking}
If a manifold of the form $\R^2 \times N^2$ equipped with product metric $g_0 \times g_N$ is a non-Bach-flat gradient Bach soliton, then it is a shrinking soliton. Furthermore, the soliton is steady if and only if it is Bach flat. 
\end{theorem}
\bp
Consider a homogeneous manifold of $\R^2 \times N^2$. Using the following equations from \cite{DK}, \cite{Ho}
\begin{equation} \label{2x2 Bii} \ba
    B_{\mu \nu} = \frac{1}{3} \nabla_{\mu} \nabla_{\nu} S_{M} - \frac{1}{3} g^M_{\mu\nu} \left[\nabla_\alpha \nabla_\alpha S_{M} - \frac{1}{2} \nabla_k\nabla_k S_{N} + \frac{1}{4} \left( \left(S_{M}\right)^2 - \left( S_{N}\right)\right)^2 \right] \text{ in }M \\
    B_{ij} = \frac{1}{3} \nabla_{i} \nabla_{j} S_{N} - \frac{1}{3} g^N_{ij} \left[\nabla_k\nabla_k S_{N} - \frac{1}{2} \nabla_\alpha \nabla_\alpha S_{M}  + \frac{1}{4} \left( \left(S_{N}\right)^2 - \left( S_{M}\right)\right)^2 \right] \text{ in }N
    \ea
\end{equation}
where $M= \R^2$, $N= N^2$, $S_M$ and $S_N$ are the respective scalar curvatures, and $g_0$ and $g_N$ are their respective metrics. Recall that homogeneous 2-manifolds have constant scalar curvature, thus we see that:
\[B_{00} = \frac{1}{12} (S_N)^2 g_{00} \quad B_{11} = \frac{1}{12} (S_N)^2 g_{11} \quad B_{22} = -\frac{1}{12} (S_N)^2 g_{22}  \quad B_{33} = -\frac{1}{12} (S_N)^2 g_{33}.\]

Since $\R^2 \times N^2$ is a gradient Bach soliton, the following system must hold.

\[ \begin{bmatrix} f_{xx}g_{00}&0&0&0 \\ 0&f_{yy}g_{11}&0&0 \\ 0&0&0&0 \\ 0&0&0&0\end{bmatrix} = \begin{bmatrix} \left(\frac{1}{24}(S_N)^2+c\right) g_{00}&0&0&0 \\ 0&\left(\frac{1}{24}(S_N)^2+c\right) g_{11} &0&0\\ 0&0&\left(\frac{-1}{24}(S_N)^2+c\right) g_{22} &0\\ 0&0&0 &\left(\frac{-1}{24}(S_N)^2+c\right) g_{33} \end{bmatrix} \]

Thus $0 =\left(\frac{-1}{24}(S_N)^2+c\right) g_{ii}$ for $i=2,3$.  Since the metric is positive definite, we know that $c= \frac{1}{24}(S_N)^2$. Thus $c \geq 0$ and the soliton must be steady or shrinking. 

The soliton is steady if and only if $S_N = 0$ which happens if and only if the manifold is Bach flat.

If the manifold is non-Bach-flat, then $c>0$ and soliton must be shrinking.
\ep

Scaling $S^2$ and $H^2$ so that $S_{S^2}= 1 = -S_{H^2}$, we see that $c= \frac{1}{24}$ and the potential function is of the form $f(x,y) = \frac{1}{24}(x+y)^2 + ax+by+d$. Again, this differs slightly from Ho because of our initial definition of a gradient Bach soliton. This confirms that the gradient solitons found by Ho are in fact the only gradient solitons on $\R^2 \times S^2$ and $\R^2 \times H^2$ up to scaling.

\begin{corollary}\label{2x2 flat}
The potential function of a steady gradient Bach soliton of the form $\R^2 \times N^2$ equipped with product metric $g_0 \times g_N$ must be linear.
\end{corollary}
\bp
Since $\R^2 \times N^2$ must be steady, we know that $f_{xx}= f_{yy} = 0$. Using calculus, it is clear that $f(x,y) = ax+by+d$. 
\ep

\begin{corollary}
The manifold $\R^2 \times \R^2$ with metric $g =g_0 \times g_N$, where $g_N$ is a flat metric, is a steady gradient Bach soliton with linear potential function. 
\end{corollary}
\bp
Consider a homogeneous manifold of $\R^2 \times \R^2$. Using (\ref{2x2 Bii}), we know that $\R^2\times \R^2$ is Bach flat. By Theorem \ref{2x2 shrinking} we know that the soliton is steady. By Corollary \ref{2x2 flat} the potential function must be linear. 
\ep

\subsection{Manifolds of the form $\R \times N^3$} 

We begin by stating and proving statements that apply to all homogeneous manifolds of the form $\R\times N^3$, then we will examine specific manifolds of this form.

A few notes before stating the theorem. We will look at a potential function $f: \R \to \R$. Since I use $x$ in later computations to mean something else, I have chosen to make $f$ a function of $r\in \R$. Furthermore, note that in this potential function $c\in \R$ is the same $c$ such that $\hess f = cg+ \frac{1}{2} B$. Thus, is we have a steady soliton, the potential function necessarily lacks that term.

\begin{lemma}\label{potential function}
A gradient Bach soliton of the form $\R\times N^3$ with metric $g=g_0 \times g_N$ has potential function of the form $f(r)= 2cr^2+ a r+b$ for $a, b\in \R$. 
\end{lemma}

\bp
Since the manifold is a soliton, we know that $\hess f = cg+ \frac{1}{2} B$. By Theorem \ref{PW Theorem} that $f$ is a function on $r\in \R$ and consequently $\tr \hess f = f''(r)$. Since the Bach tensor is trace free:
\[ \tr \hess f = \tr(cg) + \tr B \Longrightarrow f''(r)= 4c \]
Using calculus we see that this implies that $f(r)= 2cr^2+ a r+b$ for $a,b\in \R$. 
\ep

In order to examine specific manifolds, we will need the following theorem. This theorem enables us to use algebra to determine which metrics will produce solitons. 

\begin{theorem} \label{GBS}
 Consider a manifold of the form $\R \times N^3$ equipped with metric $g=g_0 \times g_N$. The manifold is a gradient Bach soliton if and only if 
\begin{equation}\label{Bii} \frac{B_{11}}{g_{11}}=\frac{ B_{22}}{g_{22}}=\frac{B_{33}}{g_{33}}= - 2c\quad\quad\text{for } c\in \R \end{equation}
\end{theorem}

\bp Consider a manifold of the form $\R \times N^3$ equipped with metric $g=g_0 \times g_N$. Suppose that this manifold is a gradient Bach soliton. Then:
\[ \hess f = cg + \frac{1}{2} B\]
where $f: \R \to \R$. Examining the components of the flow:

\[\begin{bmatrix}f''g_{00}&0&0&0 \\ 0&0&0&0 \\ 0&0&0&0 \\ 0&0&0&0\end{bmatrix} = c \begin{bmatrix}g_{00}&0&0&0 \\ 0& g_{11} &0&0 \\ 0&0& g_{22} &0 \\ 0&0&0&g_{33} \end{bmatrix}+ \frac{1}{2} \begin{bmatrix}B_{00}&0&0&0 \\ 0& B_{11} &0&0 \\ 0&0& B_{22} &0 \\ 0&0&0&B_{33} \end{bmatrix}. \]
This system yields the following equalities.
 \[ f'' g_{00} - \frac{1}{2} B_{00}= cg_{00} \quad\quad -\frac{1}{2}B_{11} = c g_{11} \quad\quad -\frac{1}{2} B_{22}= c g_{22} \quad\quad -\frac{1}{2}B_{33} = c g_{33}\]
It follows that:
\[ \frac{B_{11}}{g_{11}}=\frac{ B_{22}}{g_{22}}=\frac{B_{33}}{g_{33}}= - 2c\quad\quad\text{for } c\in \R \]
Thus the desired equality holds.\\

Further, Since $B_{00}= -2cg_{00} + 2 f''(r) g_{00} =6c g_{00}$, we see that $\frac{B_{00}}{g_{00}} = 6c$. \\ 

Suppose, on the other hand, that 
\[ \frac{B_{11}}{g_{11}}=\frac{ B_{22}}{g_{22}}=\frac{B_{33}}{g_{33}}=-2c\quad\quad\text{for } c\in \R \]
Then $-\frac{1}{2}B_{11} = c g_{11}$, $-\frac{1}{2}B_{22}= c g_{22}$, and $-\frac{1}{2}B_{33} = c g_{33}$. Taking the trace of the Bach tensor: 
\[\ba \tr B &= g^{ij}B_{ij} 
\\&= g^{00}B_{00}+ g^{11}B_{11}+g^{22}B_{22}+g^{33}B_{33} 
\\&= g^{00}B_{00}-2g^{11}c g_{11}-2g^{22}c g_{22}-2g^{33}c g_{33}
\\&= g^{00}B_{00}-6c
\ea\]
Since $B$ is trace free, we see that $B_{00} = 6c g_{00}$. By Lemma \ref{potential function} $f''(r) = 4c$, so:
\[f'' g_{00}-\frac{1}{2}B_{00}=  4c g_{00} - \frac{1}{2} (6c g_{00}) = c g_{00} \]
Thus, $ \nabla_i \nabla _j f - \frac{1}{2} B_{ij} = c g_{ij}$ for all $i, j= 0,1,2,3$, so $\hess f = cg + \frac{1}{2} B$. Therefore, $\R \times N^3$ is a gradient Bach soliton. 
\ep

From this theorem we are able to classify the resulting solitons of the form $\R \times N^3$. To do so we will need the find components of the Bach tensor using the following equation from \cite{Helliwell} and \cite{DK}.

\begin{equation} \label{1x3 Bii}
\ba 
B_{00} =& \left(-\frac{1}{12} (\Delta^{(2)} S^{(2)}) - \frac{1}{4} \left[ (|\ric|^{(2)})^2 - \frac{1}{3} (S^{(2)})^2 \right] \right) g_{00}
\\B_{jk} =& \frac{1}{2} \Delta^{(2)} \ric_{jk}^{(2)} -\frac{1}{12} \Delta^{(2)}S^{(2)} g_{jk} -\frac{1}{6} S_{;jk}^{(2)} - 2\tr^{(2)}(\ric^{(2)} \otimes \ric^{(2)})_{jk} 
\\&\;\;+ \frac{7}{6} S^{(2)}\ric_{jk}^{(2)} + \frac{3}{4}(|\ric|^{(2)})^2g_{jk} - \frac{5}{12} (S^{(2)})^2g_{jk}
\ea
\end{equation} 
where $M^{(1)} = \R$ and $M^{(2)} = N^3$. 

\begin{corollary} \label{expanding}
 If a manifold of the form $\R \times N^3$ equipped with metric $g=g_0 \times g_N$ is a non-Bach-flat gradient Bach soliton, then it is an expanding soliton. The soliton is steady if and only if it is Bach flat.
\end{corollary}

\bp
Consider a manifold of the form $\R \times N^3$ equipped with metric $g=g_0 \times g_N$. \\

From Theorem \ref{GBS} we know that:
 \[\frac{B_{11}}{g_{11}}=\frac{ B_{22}}{g_{22}}=\frac{B_{33}}{g_{33}}= - 2c \]
 
Since the Bach tensor is trace free we know that: 
\[\ba - B_{00}&= \frac{B_{11}}{g_{11}} g_{11} +\frac{ B_{22}}{g_{22}} g_{22} + \frac{B_{33}}{g_{33}} g_{33} 
\\ &= -2c (g_{11} + g_{22} + g_{33} )
\\ B_{00} &= 2c(g_{11}+ g_{22}+ g_{33})
\ea\]
Using (\ref{1x3 Bii}), since $S$ is constant:
\[B_{00} = - \frac{1}{4} \left[ (|\ric|^{(2)})^2 - \frac{1}{3} (S^{(2)})^2 \right] g_{00}\]
By Cauchy-Schwartz, we know 
\[ |\ric^{(2)}|^2 \geq \frac{\tr\left(\ric^{(2)}\right)}{ 3} = \frac{1}{3}(S^{(2)})^2,\]
and thus $B_{00} \leq 0$. Since the metric is positive definite, this implies $c \leq 0$, where $c=0$ if and only if $B_{00} =0$. By definition a soliton is expanding if $c<0$. 

If $c=0$, $B_{00} = 0$ then:
 \[ \begin{bmatrix} f_{00}&0&0&0 \\ 0&0&0&0 \\ 0&0&0&0 \\ 0&0&0&0\end{bmatrix} =  \frac{1}{2} \begin{bmatrix}0&0&0&0 \\ 0& B_{11} &0&0 \\ 0&0& B_{22} &0 \\ 0&0&0&B_{33} \end{bmatrix} \]
 Clearly, this implies that $B_{ii}=0$ for $i=1,2,3$. Thus, if the soliton is steady, the manifold is Bach flat. 
 
 If the soliton is Bach flat then $\hess f = cg$, so $0 = cg_{ii}$ for $i=1,2,3$ so $c=0$ and the soliton is steady.
\ep

\begin{remark} Recall that rescaling is a diffeomorphism of $\R$. Consequently, shrinking and expanding are diffeomorphic to one another. That is, contracting is the same as stretching after diffeomorphism. Applying this to our soliton, we see that though $\frac{\partial}{\partial t} g_{00} <0$ under the Bach flow (\cite[Proposition 2.2]{Helliwell}), $\R \times N^3$ is expanding as a soliton.
\end{remark}

In order to use this theorem to find metrics that produce solitons, we will need explicit representations of the Bach tensor. These can be found using (\ref{1x3 Bii}). The Bach tensor for solitons of the form $\R \times N^3$ where $N^3$ 3-dimensional unimodular Lie group is given in \cite{Helliwell}. For other Lie groups, one can find the necessary information using the structure constants (see \cite{Milnor} \cite{RS} \cite{IJ92}) and the equations in \cite{Helliwell} to find the necessary information for (\ref{1x3 Bii}). It should be noted that the calculations involved in finding the components of the Bach tensor are non-trivial and require the use of mathematical software.

We will begin investigating manifolds of the form $\R \times N^3$ by examining the covering spaces for the nine manifolds with compact quotient. The qualitative behavior of the compact quotients is examined in \cite{Helliwell}. The gradient solitons of the compact quotients themselves are easily classified by Corollary \ref{aof compact}. We, however, are interested in the solitons on the covering spaces themselves.

Proceeding, we will examine the 9 manifolds in \cite{Helliwell} to see if there is a metric that produces a gradient Bach solitons. The Lie groups with compact quotient are given by the unimodular, solvable Bianchi classes. That is, Bianchi classes I, II, VI$_0$, VII$_0$, VIII, and IX. There are three additional cases which are not Lie groups, but have compact quotient.  

By Theorem \ref{GBS} we need only show that a metric satisfies (\ref{Bii}). If there are no metrics that satisfy the string of equalities, then the manifold produces no solitons. The general methodology is to use the explicit representation for the Bach tensor in the above equality, then see what conditions must be placed on the metric to produce a soliton. For ease of notation in calculations, we will let:
\[ x= g_{11}, \quad y=g_{22}, \quad z=g_{33}, \quad \beta = \frac{1}{6(\det g)^2}.\] 
To clarify the consequences of each example, the metric notations will be used. These proofs heavily rely on the fact that Reimannian metrics are positive definite. That is, $g_{ii}>0$ is a strict inequality. This allows us the use the quotients in (\ref{Bii}) and to rule out potential solitons. A summary of our results is as follows. The proofs will be in subsequent sections.

\begin{theorem}
For a homogeneous manifold of type $M= \R^1 \times N^3$ equipped with the metric $g=g_0 \times g_N$ the following hold:
\bea
\item If $N^3= \R^3$, then a metric $g= g_0 \times g_N$, where $g_N$ is a flat metric, produces a gradient Bach soliton with linear potential function.
\item If $N^3 = Nil,\; Solv,\; \widehat{SL}(2, \R),\; \R \times S^2, \; \R \times H^2$ then $g$ is not a gradient Bach soliton
\item If $N^3 = E(2), H^3$, then $g$ produces a Bach soliton if and only if it is Bach flat.
\item If $N^3= S^3$, then a gradient Bach soliton is produced if and only if the metric is of the from $g_{11} = g_{22} = g_{33}$ or if it is isometric to $g_{11}= g_{22}= 4g_{33}$. These solitons are categorized in Theorems \ref{111} and \ref{114} respectively.
\ee
\end{theorem}

\subsubsection{$\R \times \R^3$}

\begin{prop}
The manifold $\R \times \R^3$ with metric $g =g_0 \times g_N$, where $g_N$ is a flat metric, is a gradient Bach soliton with potential function $f(r) = ar+b$ or some $a \in \R$.
\end{prop}

\bp
We know from (\ref{1x3 Bii}) that $B_{ii} = 0$ for $i=0,1,2,3$. By Corollary \ref{expanding} we know that the soliton is steady, so $c=0$. So by Lemma \ref{potential function} $f(r) = ar+b$ for $a, b\in \R$. 
\ep

\subsubsection{$\R \times Nil$}

We know from \cite{Helliwell}
    \[\ba B_{00}= -\beta (g_{00})^3 (g_{11})^4 \quad & \quad B_{11} = -5\beta(g_{00})^2 (g_{11})^5  
    \\ B_{22} = 3\beta(g_{00})^2 (g_{11})^4 g_{22} \quad & \quad B_{33} = 3\beta(g_{00})^2 (g_{11})^4 g_{33}. \ea\]

\begin{prop}
The manifold $\R \times Nil$ with metric $g= g_0 \times g_{Nil}$ is not a gradient Bach soliton.
\end{prop}
\bp
Proceeding by contradiction, suppose $\R \times Nil$ with metric $g= g_0 \times g_{Nil}$ is a gradient Bach soliton. Then using (\ref{Bii}) we see that:
\[ \frac{B_{11}}{g_{11}} = \frac{B_{22}}{g_{22}} \quad\Longrightarrow \quad -5\beta (g_{00})^2(g_{11})^4 = 3 \beta (g_{00})^2(g_{11})^4.\]
However, this implies that $-5 = 3$. Thus $\R \times Nil$ is not a gradient Bach soliton. 
\ep
    
\subsubsection{$\R \times Solv$}

We know from \cite{Helliwell}
    \[\ba  B_{00} = - \beta p(g_{11}, g_{22}) (g_{00})^3 \quad & \quad  B_{11} = -\beta q(g_{11}, g_{22}) (g_{00})^2g_{11}  
    \\ B_{22} = -\beta q(g_{22}, g_{11}) (g_{00})^2g_{22} \quad & \quad B_{33}= 3\beta p(g_{11}, g_{22}) (g_{00})^2g_{33}\ea\]
    where
    \[ p(x,y) = x^4 + x^3y + xy^3 +y^4 \quad\quad q(x,y) = 5x^4 + 3x^3y - xy^3 -3y^4.\]

\begin{prop}
The manifold $\R \times Solv$ with metric $g= g_0\times g_{Solv}$ is not a gradient Bach soliton.
\end{prop}
\bp
Proceeding by contradiction, suppose $\R \times Solv$ with metric $g= g_0\times g_{Solv}$ is a gradient Bach soliton. Using (\ref{Bii}) we see that: 
\[ \frac{B_{11}}{g_{11}} = \frac{B_{33}}{g_{33}} \quad\Longrightarrow \quad -\beta\;q(g_{11}, g_{22}) (g_{00})^2 = 3 \beta\; p (g_{11}, g_{22}) (g_{00})^2 \]
Letting $x= g_{11}$ and $y= g_{22}$:
\[ \bas -q(x,y) &= 3p(x,y) 
\\ -5x^4 - 3x^3y + xy^3+ 3y^4 &= 3x^4+ 3x^3y +3xy^3 + 3y^4
\\ -2x(4x^3+6x^2y + y^3)&=0
\eas\]
Then either $x=0$ or $4x^3+6x^2y + y^3 =0$. The first statement is not possible because the metric is positive definite. The latter statement holds if and only if $x=y=0$ forcing either $g_{11}=0$ or $g_{11}= g_{22}= 0$, contradicting positive definiteness. Thus $\R \times Solv$ is not a gradient Bach soliton. 
\ep

\subsubsection{$\R \times \widehat{SL} (2, \R)$}

We know from \cite{Helliwell}
    \[\ba  B_{00} = -\beta p(-g_{11}, g_{22}, g_{33})(g_{00})^3 \quad & \quad  B_{11} = -\beta q(-g_{11}, g_{22}, g_{33}) (g_{00})^2g_{11}  
    \\B_{22} = -\beta q(g_{22}, -g_{11}, g_{33}) (g_{00})^2g_{22} \quad & \quad B_{33}= -\beta q(g_{33}, -g_{11}, g_{22}) (g_{00})^2g_{33} \ea\]
    where
    \[\bas p(x,y,z) &= x^4- x^3(y+z)+ x^2yz+ x(-y^3+y^2z+ yz^2-z^3) + y^4-y^3z-yz^3+z^4 
    \\ q(x,y,z)&= 5x^4- 3x^3(y+z)+ x^2yz+ x(y^3-y^2z- yz^2+z^3) -3y^4+ 3y^3z+3yz^3-3z^4. \eas\]

\begin{prop}
The manifold $\R \times \widehat{SL}(2, \R)$ with metric $g= g_0\times g_{\widehat{SL}(2, \R)}$ cannot be a gradient Bach soliton.
\end{prop}
\bp
Proceeding by contradiction, suppose $\R \times \widehat{SL}(2,\R)$ with metric $g= g_0\times g_{\widehat{SL}(2, \R)}$ is a gradient Bach soliton. Using (\ref{Bii}) we see that:

\[\bas \frac{B_{22}}{g_{22}} &= \frac{B_{33}}{g_{33}}
\\q(y,-x, z) &=q(z,-x, y)
\\ \bar
5y^4+3xy^3-3y^3z- xy^2z-x^3y-x^2yz \\
       + xyz^2+yz^3-3x^4- 3x^3z- 3xz^3-3z^4 
\end{array}
 \quad &=\quad
 \bal
       5z^4+3xz^3-3yz^3- xyz^2-x^3z-x^2yz \\
      + xzy^2+y^3z-3x^4- 3x^3y- 3xy^3-3y^4 
 \end{array}
 \\ \bar
  2 (y - z) (x^3 + 3 x y^2 + 2 x y z + 3 x z^2 
\\+ 4 y^3 + 2 y^2 z + 2 y z^2 + 4 z^3)
\end{array}
\quad &=\quad 0
 \eas\]
 The only potential real solution is that $y=z$. As above, because the metric is positive definite, the last term in the product is nonzero. Examining the consequences of this using the other equations in (\ref{Bii}) we see that the following must hold.

\[\bas \frac{B_{11}}{g_{11}} &= \frac{B_{22}}{g_{22}} 
\\ q(-x, y, z) &=q(y,-x, z)
\\ \bar
 5x^4+3x^3y+3x^3z + x^2yz - xy^3+xy^2z 
\\ +xyz^2-xz^3 -3y^4+ 3y^3z+3yz^3-3z^4 
\end{array}
 \quad &=\quad
 \bal
       5y^4+3xy^3-3y^3z- xy^2z-x^3y-x^2yz \\
      + xyz^2+yz^3-3x^4- 3x^3z- 3xz^3-3z^4 
 \end{array}
 \\ \bar
8x^4+ 4x^3y +6x^3z +  2x^2yz - 4xy^3
\\+2xy^2z  +2xz^3 -8y^4+ 6y^3z+2yz^3 
\end{array}
 \quad &=\quad 0
 \eas\]
However, if $y=z$ then: 
 \[\bas
 \bar
 8x^4+ 4x^3y +6x^3z +  2x^2yz - 4xy^3
\\+2xy^2z  +2xz^3 -8y^4+ 6y^3z+2yz^3 
\end{array}
 \quad &=\quad 
 \bal 8x^4 + 4x^3y+6x^3y +  2x^2y^2 - 4xy^3
\\ +2xy^3 +2xy^3 -8y^4+ 6y^4+2y^4
\end{array}
\\  \; & = \quad 8x^4+ 10x^3y + 2x^2y^2
\\\;&\neq 0
 \eas\]
 
Therefore if $y=z$, then $B_{11}\;/ \;g_{11} \neq B_{22}\;/ \;g_{22}$. Thus $y\neq z$. Therefore, $\R \times \widehat{SL}(2,\R)$ is not a gradient Bach soliton.
\ep

\subsubsection{$\R \times (\R \times S^2)$} 

\begin{prop}
There are no gradient Bach solitons on $\R \times (\R \times S^2)$ with metric $g= g_0 \times (g_\R \times g_{S^2})$.
\end{prop}

\bp 
Consider the manifold $\R\times (\R \times S^2)$ with metric $g= g_0 \times (g_\R \times g_{S^2})$. Rescaling the sphere to have scalar curvature $S_{S^2} = 1$, from Theorem \ref{2x2 shrinking} we know:
\[B_{00} =\frac{1}{12}g_{00}\quad\quad B_{11} =\frac{1}{12}g_{11}\quad\quad B_{22} =-\frac{1}{12}g_{22}\quad\quad B_{33} =-\frac{1}{12}g_{33}.  \]
This contradicts Theorem \ref{GBS}. Therefore, there are no gradient Bach solitons on $\R \times (\R \times S^2)$ with potential function on $\R$. 

\ep

\subsubsection{$\R \times (\R\times H^2)$}
\begin{prop} 
There are no gradient Bach solitons on $\R \times (\R \times H^2)$ with metric $g= g_0 \times (g_\R \times g_{H^2})$.
\end{prop}

\bp 
Rescaling the $H^2$ to have scalar curvature $S_{H^2} = -1$, from Theorem \ref{2x2 shrinking} we know:
\[B_{00} =\frac{1}{12}g_{00}\quad\quad B_{11} =\frac{1}{12}g_{11}\quad\quad B_{22} =-\frac{1}{12}g_{22}\quad\quad B_{33} =-\frac{1}{12}g_{33},\]
and thus the proof follows exactly as in the proof for $\R \times \R \times S^2$ above.
\ep

\subsubsection{$\R \times E(2)$} 

We know from \cite{Helliwell}
    \[\ba  B_{00} = -\beta p(-g_{11}, g_{22}) (g_{00})^3 \quad & \quad  B_{11} = -\beta q(-g_{11}, g_{22}) (g_{00})^2g_{11}  
    \\ B_{22} = -\beta q(g_{22}, -g_{11}) (g_{00})^2g_{22} \quad& \quad B_{33}= 3\beta p(-g_{11}, g_{22}) (g_{00})^2g_{33} \ea\]
    where $p(x,y)$ and $q(x,y)$ are as above. 
    
\begin{prop}
The manifold $\R \times E(2)$ with metric $g= g_0 \times g_{E(2)}$ is a gradient Bach soliton if and only if it is Bach flat.
\end{prop}
\bp
Consider the manifold $\R \times E(2)$ with metric $g= g_0 \times g_{E(2)}$. Using (\ref{Bii}) we see that:
\[\bas \frac{B_{11}}{g_{11}} &= \frac{B_{22}}{g_{22}} 
\\q(-x, y) &=q(y,-x)
\\5x^4 -3x^3y +xy^3-3y^4 &= 5y^4- 3y^3x+ yx^3-3x^4
\\(x-y)(x+y)(2x^2-xy+2y^2) &= 0
\eas\]
The only two real, nonzero solutions are that $x=y$ or $x=-y$. Since our metric is positive definite $x\neq -y$. Thus $x=y$ is the only candidate. Proceeding, we will see that the equalities from (\ref{Bii}) are satisfied if and only if $x=y$. 

\[ \bas\frac{B_{11}}{g_{11}} &= \frac{B_{33}}{g_{33}} 
\\-q(-x, y) &= 3p(-x,y) 
\\-5x^4 + 3x^3y - xy^3 + 3y^4 &= 3x^4 - 3x^3y -3xy^3 +3y^4
\\-2x( 4x^3- 3x^2y- y^3) &=0
\eas\]
Since $x\neq 0$, $4x^3 - 3x^2y - y^3 = 0$. We see that $x=y$ holds. 
\[ \bas\frac{B_{22}}{g_{22}} &= \frac{B_{33}}{g_{33}} 
\\-q(y, -x) &= 3p(-x,y) 
\\-5y^4 + 3xy^3-x^3y+3x^4 &= 3x^4 - 3x^3y -3xy^3 +3y^4 
\\-2y( 4y^3- 3xy^2- x^3) &=0
\eas\]
Since $y\neq 0$, $4y^3- 3xy^2- 2x^3 =0$. Again, we see that $x=y$ holds. \\

Thus, $g_{11} = g_{22}$. This condition is equivalent to being Bach flat by the following lemma. Therefore, by Theorem \ref{GBS} and Lemma \ref{E2 flat}, $\R \times E(2)$ is a gradient Bach soliton if and only if it is Bach flat. 
\ep

\begin{lemma} \label{E2 flat}
The manifold $\R \times E(2)$ with metric $g= g_0 \times g_{E(2)}$ is Bach flat if and only if $g_{11} = g_{22}$.
\end{lemma}

\bp
Factoring the components of the Bach tensor for $\R \times E(2)$: 
\[\bas  B_{00} &=    -\beta \; (g_{11}-g_{22})^2\left( (g_{11})^2+ g_{11}g_{22} +(g_{22})^2\right)(g_{00})^3 
\\ B_{11} &=  -\beta \; (g_{11}-g_{22})\left( 5(g_{11})^3+ 2(g_{11})^2(g_{22}) + 2 (g_{11})(g_{22})^2+3(g_{22})^3 \right)(g_{00})^2 g_{11} 
\\ B_{22} &= - \beta \; (g_{22}-g_{11})\left( 3(g_{11})^3+ 2(g_{11})^2(g_{22}) + 2 (g_{11})(g_{22})^2+3(g_{22})^3 \right)(g_{00})^2g_{22} 
\\B_{33} & =   3\beta \; (g_{11}-g_{22})^2\left( (g_{11})^2+ g_{11}g_{22} +(g_{22})^2\right)(g_{00})^2 g_{11} 
\eas\]
Since our metric is positive definite $B_{ii} =0$ if and only if $g_{11} - g_{22} = 0 $ if and only if $g_{11} = g_{22}$.
\ep

\subsubsection{$\R \times H^3$} 

\begin{prop}
The manifold $\R \times H^3$ with metric $g= g_0 \times g_{H^3}$ is the trivial gradient Bach soliton. That is, $\R \times H^3$ is a Bach soliton if and only if it is Bach-flat.
\end{prop}

\bp
Following the explanation from \cite{Helliwell}, we know that $H^3$ is a one parameter family of homogeneous metrics. Consequently all metrics are Einstein since they are scalar multiples of the standard metric. Thus, as Helliwell concludes, the flat metric remains flat in the Bach flow. Therefore, the Bach flat metric produces a gradient soliton.
\ep

\subsubsection{$\R \times S^3$}

Before delving into this case, it is important that the reader note that I am $S^3$ to be synonymous with $SU(2)$. That is, the manifold does NOT necessarily have the round metric, but rather has any left invariant metric on Lie group $SU(2)$. My choice to call this $S^3$ was motivated by wanting to maintain consistency between the cases presented by Helliwell in \cite{Helliwell} and this paper.

We know from \cite{Helliwell}
    \[\ba  B_{00} = -\beta\; p(g_{11}, g_{22}, g_{33})(g_{00})^3 \quad & \quad B_{11} = -\beta\; q(g_{11}, g_{22}, g_{33}) (g_{00})^2g_{11}  
    \\B_{22} = -\beta \;q(g_{22}, g_{33}, g_{11}) (g_{00})^2g_{22} \quad & \quad B_{33}= -\beta\; q(g_{33}, g_{11}, g_{22}) (g_{00})^2g_{33} \ea\]
    where
    \[\bas p(x,y,z) &= x^4- x^3(y+z)+ x^2yz+ x(-y^3+y^2z+ yz^2-z^3) + y^4-y^3z-yz^3+z^4 
    \\ q(x,y,z)&= 5x^4- 3x^3(y+z)+ x^2yz+ x(y^3-y^2z- yz^2+z^3) -3y^4+ 3y^3z+3yz^3-3z^4 \eas\]

\begin{prop} \label{S3}
The manifold $\R \times S^3$ with metric $g= g_0 \times g_{SU(2)}$ is a gradient Bach soliton if and only if our metric is $g_{11} = g_{22} = g_{33}$ or if it is isometric to $g_{11} = g_{22} = 4g_{33}$. 
\end{prop}

\bp
Proceeding, consider $\R \times S^3$ with metric $g= g_0 \times g_{SU(2)}$. We will show that the (\ref{Bii}) holds if and only if $x=y=z$, $x=y=4z$, $x=4y=z$, or $4x= y=z$.\\

We will first consider that case where $x=y=z$: 
\[  \frac{B_{11}}{g_{11}}   = \frac{B_{22}}{g_{22}}= \frac{B_{33}}{g_{33}}= -\beta\; q(g_{11}, g_{11}, g_{11}) (g_{00})^2 \]
This clearly satisfies (\ref{Bii}). 

Proceeding to examine the equalities in general we see that:

 \begin{equation}\label{S3 11 22} 
 \bas 
 \frac{B_{11}}{g_{11}} &= \frac{B_{22}}{g_{22}}&
\\q(x, y, z) &=q(y,z, x) &
\\ \bar
 5x^4- 3x^3y-3x^3z+ x^2yz+ xy^3-xy^2z
\\ - xyz^2+xz^3 -3y^4+ 3y^3z+3yz^3-3z^4 
\end{array}
 \quad &=
 \bal
      5y^4- 3y^3z-3xy^3+ xy^2z+ yz^3-xyz^2
    \\ - x^2yz+x^3y -3z^4+ 3xz^3+3x^3z-3x^4 
 \end{array} \\ 
 \bar
 2 (x - y) (4 x^3 + 2 x^2 y - 3 x^2 z + 2 x y^2
\\- 2 x y z + 4 y^3 - 3 y^2 z - z^3)
\end{array}
\quad &=\quad 0
\eas \end{equation}

  \begin{equation}\label{S3 11 33} 
\bas 
\frac{B_{11}}{g_{11}} &= \frac{B_{33}}{g_{33}}& 
\\q(x, y, z) &=q(y,z, x) &
\\ \bal
 5x^4- 3x^3y-3x^3z+ x^2yz+ xy^3-xy^2z 
\\ - xyz^2+xz^3 -3y^4+ 3y^3z+3yz^3-3z^4 
\end{array}
 \quad &= 
 \bal
       5z^4- 3xz^3 - 3yz^3+ xyz^2+ x^3z-x^2yz 
\\ - xy^2z+y^3z -3x^4+ 3x^3y+3xy^3-3y^4
 \end{array}
\\ \bar
 2 (x - z) (4 x^3 - 3 x^2 y + 2 x^2 z- 2 x y z 
\\+ 2 x z^2 - y^3 - 3 y z^2 + 4 z^3) 
\end{array}
 \quad &=\quad 0
 \eas 
 \end{equation}

 \begin{equation} \label{S3 22 33}
 \bas \frac{B_{22}}{g_{22}} &= \frac{B_{33}}{g_{33}} 
\\q(y,z,x) &=q(y,z, x)
\\ \bar
5y^4- 3y^3z-3xy^3+ xy^2z+ yz^3-xyz^2 
    \\ - x^2yz+x^3y -3z^4+ 3xz^3+3x^3z-3x^4
\end{array}
 \quad &=
 \bal
       5z^4- 3xz^3 - 3yz^3+ xyz^2+ x^3z-x^2yz
\\ - xy^2z+y^3z -3x^4+ 3x^3y+3xy^3-3y^4
 \end{array} 
\\ \bar
 -2 (y - z) (x^3 + 3 x y^2 + 2 x y z+ 3 x z^2 
\\- 4 y^3 - 2 y^2 z - 2 y z^2 - 4 z^3)
\end{array}
 \quad &=\quad 0 
 \eas \end{equation}
 
\begin{case} Suppose that $x=y$. Then (\ref{S3 11 22}) is satisfied. Moreover this means that in order for (\ref{S3 11 33}) to be satisfied: 
 \[\bas 0&=  4 x^3 - 3 x^3 + 2 x^2 z- 2 x^2 z + 2 x z^2 - x^3 - 3 x z^2 + 4 z^3 
 \\&= z^2(4z-x)  \eas\]
 Consequently $x=4z$. We see that this equality not only holds in $\ref{S3 22 33}$, but is forced:
 \[ \bas 0&=  x^3 + 3 x^3 + 2 x^2 z+ 3 x z^2 - 4 x^3 - 2 x^2 z - 2 x z^2 - 4 z^3
 \\&= z^2(x-4z) \eas \]
 Thus $x=y=4z$ maintains all three equalities.
 \end{case} 
 
\begin{case} Suppose that $x=z$. Then (\ref{S3 11 33}) is satisfied. Moreover this means that in order for (\ref{S3 11 22}) to be satisfied: 
 \[\bas 0&=  4 x^3 + 2 x^2 y - 3 x^3 + 2 x y^2- 2 x^2 y + 4 y^3 - 3 y^2 x - x^3
 \\&= y^2(4y-x)  \eas\]
 Consequently $x=4y$. We see that this equality not only holds in (\ref{S3 22 33}), but is forced:
 \[ \bas 0&=  x^3 + 3 x y^2 + 2 x^2 y + 3 x^3 - 4 y^3 - 2x y^2  - 2 x^2 y - 4 x^3
 \\&= y^2(x-4y) \eas \]
Thus $x=4y=z$ maintains all three equalities.
\end{case}

\begin{case} Suppose that $y=z$. Then (\ref{S3 22 33}) is satisfied. Moreover this means that in order for (\ref{S3 11 22}) to be satisfied: 
 \[\bas 0&=  4 x^3 + 2 x^2 y - 3 x^2 y + 2 x y^2- 2 x y^2+ 4 y^3 - 3 y^3 - y^3
 \\&= x^2(4x-y)  \eas\]
 Consequently $4x=y$. We see that this equality not only holds in (\ref{S3 22 33}), but is forced:
 \[ \bas 0&= 4 x^3 - 3 x^2 y + 2 x^2 y- 2 x y^2 + 2 x y^2 - y^3 - 3 y^3 + 4 y^3
 \\&= x^2(4x-y) \eas \]
 Thus $4x=y=z$ maintains all three equalities.
 \end{case}
 
\begin{case} Suppose that $x\neq y$, $x\neq z$, $y\neq z$. Then only other permissible metric would need to satisfy the system of equations:
\[\begin{cases}
4x^3 + 2x^2y -3x^2z + 2xy^2 - 2xyz + 4y^3 - 3y^2z - z^3 =0\\
4x^3 - 3x^2y + 2x^2 z - 2xyz + 2xz^2 - y^3 - 3yz^2 + 4z^3 = 0 \\
x^3 + 3 x y^2 + 2 x y z+ 3 x z^2- 4 y^3 - 2 y^2 z - 2 y z^2 - 4 z^3 = 0
\end{cases}\] 
Subtracting the first equation from the second yields:
\[\bas 5x^2y - 5x^2z +2xy^2 -2xz^2 +5y^3 -3y^2z +3yz^2 -5z^3&=0 
\\ (y - z) (5 x^2 + 2 x y + 2 x z + 5 y^2 + 2 y z + 5 z^2) &= 0\eas\]

Thus $y=z$ contradicting the original assertion. Moreover, the metric is positive definite. Thus, this case yields no potential metrics.
\end{case}

Therefore, by Theorem \ref{GBS}, $\R \times S^3$ is a Bach soliton if and only if $g_{11}= g_{22} = g_{33}$, $g_{11} = g_{22} = 4g_{33}$, $g_{11}= 4g_{22} = g_{33}$, or $4g_{11} = g_{22} = g_{33}$. 
\ep
 
 \begin{theorem}\label{111}
 If $g_{11} = g_{22} = g_{33}$ then the soliton produced by $\R\times S^3$ is Bach flat and steady.
 \end{theorem}
 \bp
Suppose $g_{11} = g_{22} = g_{33}$. We know by Theorem \ref{S3} that this is the metric of a soliton on $\R\times S^3$. Then:
\[  \frac{B_{11}}{g_{11}}   = \frac{B_{22}}{g_{22}}= \frac{B_{33}}{g_{33}}= -\beta\; q(g_{11}, g_{11}, g_{11}) (g_{00})^2 = -\beta (0) (g_{00})^2 = 0 \]
Thus $c=0$, so the soliton is steady.

Moreover, since 
    \[\bas p(x,x,x) &= x^4- x^3(2x)+ x^4+ x(-x^3+x^3+ x^3-x^3) + x^4-x^4-x^4+x^4 = 0  
    \\ q(x,x,x)&= 5x^4- 3x^3(2x)+ x^4+ x(x^3-x^3- x^3+x^3) -3x^4+ 3x^4+3x^4-3x^4=0 \eas\]
    We know that $B_{ii} = 0$ for all $i=0,1,2,3$. Therefore the metric is Bach flat.
 \ep
 
 \begin{remark}
Note that in the previous proof, one could have referenced Corollary \ref{expanding} instead of calculating the Bach tensor. The calculation was included to demonstrate an alternate method in that works when you know the components of the Bach tensor.
 \end{remark}
 
  \begin{theorem}\label{114}
 If $g_{11} = g_{22} = 4g_{33}$ then the soliton produced by $\R\times S^3$ is expanding and immortal.
 \end{theorem}
 \bp
Without loss of generality, suppose $g_{11} \leq g_{22} \leq g_{33}$. Consider $g_{11} = g_{22} = 4g_{33}$. We know by Theorem \ref{S3} that this is the metric of a soliton on $\R\times S^3$. Then:
\[  \frac{B_{11}}{g_{11}}   = \frac{B_{22}}{g_{22}}= \frac{B_{33}}{g_{33}}= -\beta\; q(g_{11}, g_{11}, 4g_{11}) (g_{00})^2 = -2 c\]
Observe that:
\[\bas q\left(x,x,\frac{1}{4}x\right) &= 5x^4- 3x^3\left(\frac{5}{4}x\right)+ \frac{1}{4}x^4+ x\left(x^3-\frac{1}{4}x^3- \frac{1}{16}x^3+\frac{1}{64}x^3\right) -3x^4+ \frac{3}{4}x^4+\frac{3}{64}x^4-\frac{3}{256}x^4
\\&= x^4 \left( 5- \frac{15}{4} + \frac{1}{4} + 1- \frac{1}{4} - \frac{1}{16} + \frac{1}{64} - 3 + \frac{3}{4} + \frac{3}{64} - \frac{3}{256} \right)
\\&= -\frac{3}{256} x^4
\eas\]
Thus $\beta \frac{3}{256} (g_{11})^4 (g_{00})^2 >0$. Since
\[ -2c =\beta \frac{3}{256} (g_{11})^4 (g_{00})^2\]
we see that $c<0$. Recall the soliton is of the form $\hess f - \frac{1}{2} B = cg$. Thus, the soliton with the given metric is expanding.

Using Theorem \ref{solutions}. The Bach tensor is conformally invariant of weight $w=-2$, so $\tau_t = \sqrt{1-4ct}$. Since $c<0$, we see that $\tau_t$ is defined for $t\in \left(\frac{1}{4c}, \infty\right)$. Thus the soliton is immortal.
 \ep 
 
\begin{remark}
This result aligns with the analysis of the Bach flow of $\R \times S^3$ in \cite{Helliwell}.
\end{remark}

\subsubsection*{Acknowledgements}
The author would like to thank Professor Dylan Helliwell of Seattle University and Professor Peter Petersen of UCLA for their interest and helpful discussions when writing this paper. Thank you also to my thesis advisor, Professor William Wylie of Syracuse University, for his guidance, support, and insight into this topic.

\bibliographystyle{alpha}
\bibliography{references}

\end{document}